\documentclass[12pt,a4paper]{article}
\usepackage{moreverb,graphicx,url}
\begin{document}
\newcommand{\hvec}[1]{\mbox{\boldmath$#1$}}
\newcommand{\re}{\mathcal{R}}
\newcommand{\ep}{\delta_{s}}
\newcommand{\pe}{\mathcal{P}}
\newcommand{\ham}{\mathcal{H}}
\newcommand{\bond}{\mathcal{B}}

\title{An Accurate Lubrication Model of Contaminated Coating Flows}
\author{A.J.~Roberts \& M.E. ~Simpson\footnote{{\protect
\url{mailto:aroberts@usq.edu.au,simpsonm@usq.edu.au}}}\\
{\small\em Department of Mathematics and Computing, University of} \\ 
{\small\em Southern Queensland, Toowoomba, Australia}}
\maketitle

\begin{abstract}
    The levelling of short-wave irregularities on a thin film of fluid 
    is primarily due to the action of surface tension. Surface tension 
    gradients are often created by a number of different factor including 
    evaporation, thermal gradients or deposition of surfactants.
    Lubrication theory, which ignores inertia terms in favour of viscous 
    terms, produces a system of two nonlinear \textsc{pde}'s for the unsteady 
    flow of a thin viscous Newtonian fluid containing an insoluble surfactant.
    A complex model, which systematically includes all relevant effects to 
    these \textsc{pde}'s, is developed using centre manifold techniques.
    The benefits of using these techniques to develop accurate models are in
    their application. Subtle variations in parameters or assumptions are
    able to be catered for by including or deleting the relevant terms rather
    than having to redevelop these models.  
    Computed solutions of both models using the same numerical process are 
    compared. Numerical simulations also demonstrate the long-term 
    stabilisation of corrugations by induced surfactant variations.

    \paragraph{Keywords:} surfactant, centre manifold, low-dimensional 
    modelling.
\end{abstract}

\tableofcontents

\section{Introduction}
Thin film flows are common in a large number of industrial and biological 
flows. Industrial applications include liquid agrochemicals, production of 
photographic film, lubricants, adhesives, dyes and surfactants. Biological 
flows include thin liquid films on the cornea of the eye and on the linings 
of the lungs. The development of accurate models is therefore essential for a 
proper understanding of these flows.

The development of models for the evolution of a thin clean film on arbitrarily
curved substrates has been well documented over recent years.
Levich~\cite{levich62} developed a model for the motion of thin and
wide fluid films induced by a surface tension variation. Inconsistencies in the
assumptions in Levich's solution were corrected by Yih~\cite{yih68} in 1968.
Kennings~\cite{kenning68} provided a valuable contribution to the qualitative
effects on interfacial motion by surface tension gradients. Ahmad \& 
Hansen~\cite{ahmad72} considered the spreading of a monolayer over a thin 
liquid film and argued that the distance spread in time $t$ is 
$x^{2}=(2H/\mu)\pi_{0}t$ where $H$ is the liquid-film thickness,
$\mu$ the coefficient of viscosity of the liquid underlying the monolayer and
$\pi_{0}$ is the spreading pressure of the lens generating the monolayer. 
This relation was verified in experiments by Hussain, Fatima \& 
Ahmad~\cite{hussain75} in 1975. 

DiPietro, Huh \& Cox~\cite{dipietro78}, DiPietro \& Cox~\cite{dipietro80} and 
Foda \& Cox~\cite{foda80} analysed the spreading of one liquid on the surface 
of a deep fluid. They provide an extensive overview of the interfacial 
dynamics associated with the spreading of a contaminant.

These papers used methods to develop models with specific parameter values. If
the physics of the problem required the reordering of these parameterised 
physical processes, then, in traditional approaches all the modelling needs to
be performed again. However, modern dynamical systems theory provides 
systematic methods to derive comprehensive and flexible low-dimensional models
of spatio-temporal evolution. Centre manifold theory is one such method used
successfully in deriving flexible and accurate models as shown by
Roberts in 1996 and 1997~\cite{roberts96,roberts97}. 
Roy {\it et al}~\cite{roy97} also used these techniques to demonstrate the 
importance of higher order terms in conserving mass.

This work continues on from that analysis by now developing thin fluid
film models including effects due to the contamination of the surface with an 
insoluble surfactant. In addition the effects of inertia and van der Waals
forces are included for completeness. The relevant equations, boundary 
conditons and constitutive equations are analysed in Section 3 using centre
manifold theory in order to generate the low-dimensional model for the 
dynamics. In Section 4 these equations are solved using the computer algebra 
package {\sc reduce} to compute evolution equations for the given
problem to any order required. The evolution equations to low order are
\begin{eqnarray}
    \frac{\partial\eta}{\partial t} & = & -\frac{1}{2}\frac{\partial}
    {\partial x}\left(\eta^{2}\gamma_{x}\right)-\frac{1}{3}\frac{\partial}
    {\partial x}\left(\eta^{3}\left(\gamma\eta_{xx}\right)_{x}\right)
    \nonumber \\
    & & {}-\frac{1}{3}\frac{\partial}{\partial x}\left(\eta^{3}\right)\bond
    \sin{\theta}+\frac{1}{3}\frac{\partial}{\partial x}\left(\eta^{3}\eta_{x}
    \right)\bond\cos{\theta}-\frac{\partial}{\partial x}\left(\frac{\eta_{x}}
    {\eta}\right)\ham\nonumber \\
    & & {}+\mathcal{O}(\partial^{6}_{x}+\bond^{2}+\ham^{2})\quad\mbox{and}
    \label{introeta} \\
    \frac{\partial\Gamma}{\partial t} & = & -\frac{\partial}{\partial x}
    \left(\Gamma\eta\gamma_{x}\right)-\frac{1}{2}\frac{\partial}{\partial x}
    \left(\Gamma\eta^{2}\left(\gamma\eta_{xx}\right)_{x}\right)\nonumber \\
    && {}-\frac{1}{2}\frac{\partial}{\partial x}\left(\Gamma\eta^{2}\right)
    \bond\sin{\theta}+\frac{1}{2}\frac{\partial}{\partial x}\left(\Gamma
    \eta^{2}\eta_{x}\right)\bond\cos{\theta}-\frac{3}{2}\frac{\partial}
    {\partial x}\left(\frac{\Gamma\eta_{x}}{\eta^{2}}\right)\ham\nonumber \\
    & & {}+\frac{\ep}{\sqrt{1+\eta_{x}^{2}}}\frac{\partial}{\partial x}\left(
    \frac{\Gamma_{x}}{1+\eta_{x}^{2}}\right)\nonumber \\
    && {}+\mathcal{O}(\partial^{6}_{x}+\bond^{2}+\ham^{2})
    \label{introgamma}
\end{eqnarray}
where $\bond$, $\ham$ and $\ep$ is a Bond number, a Hamaker constant and an
inverse P\'{e}cl\'{e}t number respectively. We have reproduced the evolution
models of Gaver and Grotberg~\cite{gaver90}, de Wit~\cite{dewit94} and others
except for the additional term in the evolution equation of the contaminant
which reflects the importance of gradient effects. However, in Section 4 we
show that these models do not capture all the terms required for accuracy. The
relevance of these terms will become apparent in the physics of the problem.

The comprehensive structurally stable model~\ref{introeta}--\ref{introgamma} 
which is
is compared numerically in Section 7 to a number of similar models developed 
recently~\cite{dewit94,gaver90,schwartz95} using traditional lubrication 
methods. These simulations show that our model better captures the effects of 
steep gradients. This is expected as only our model includes this term which 
demonstrates the flexibility of our approach, which rests on the Approximation 
Theorem~\cite{carr81}, when analysing flows with parameters in a different 
physical region. Inappropriate terms are simply deleted from the comprehensive 
model rather than having to redevelop a new model using the new assumptions. 

Finally numerical simulations of examples shown in Section 6 confirm the 
long-lasting corrugations predicted by the linear analysis in Section 5. These
corrugations only decay by the very slow diffusion of surfactant.

\section{Mathematics models the physical processes}
Consider a thin film of fluid of varying thickness $\eta(x,t)$ lying on a
flat substrate. An orthogonal coordinate system $(x,y)$
is used where $x$ measures distance
along the substrate located at $y=0$. The surface
of the fluid is thus described by $y=\eta (x,t)$.
The incompressible Newtonian fluid of viscosity $\mu$ and density $\rho$,
undergoes slow creeping flow. The dynamics of the fluid flow are determined 
by pressure gradients caused by surface tension forces which in turn are 
affected by the surfactant concentration. The equations of motion to be solved
are the continuity and Navier-Stokes equations together with boundary 
conditions. For convenience we non-dimensionalise by scaling variables with 
respect to: a typical film thickness $H$, a reference value of the surfactant 
concentration $\Gamma_{0}$ and surface tension $\gamma_{0}$ at 
$\Gamma=\Gamma_{0}$, the reference time $\mu H/\gamma_{0}$, the reference 
velocity $U=\gamma_{0}/\mu$, and the reference pressure $\gamma_{0}/H$.
The equations of motion then are written
\begin{equation}
\hvec{\nabla}\cdot\hvec{q}=0\qquad\mbox{and}
\label{continuity}
\end{equation}
\begin{equation}
\re\left[\frac{\partial\hvec{q}}{\partial t}+\hvec{q}\cdot\hvec{\nabla}
\hvec{q}\right]=-\hvec{\nabla}(p+W)+\nabla^{2}\hvec{q}+\bond\hvec{g}\, ,
\label{navier}
\end{equation}
where $\re=\gamma_{0}\rho H/\mu^{2}$ is a Reynolds number, $\bond=\rho gH^{2}/
\gamma_{0}$ is a Bond number, $W=\ham/\eta^3$ is the simple model of the van 
der Waals force used by de Wit {\it et al}~\cite{dewit94} and described in
detail by Maldarelli {\it et al}~\cite{maldarelli80} in which $\ham=H_{a}\rho/
H\mu^{2}$ is a nondimensional Hamaker constant where $H_{a}=10^{-12}$ erg is a
typical value of the dimensional Hamaker constant, $\hvec{q}(x,y,t)$ is the 
fluid velocity, $p(x,y,t)$ is the pressure field and $\hvec{g}$ is the 
direction of gravitational normal force at an 
angle $\theta$ to the substrate ($\theta=\pi/2$ is draining flow and 
$\theta=\pi$ generally leads to dripping). These equations are to be solved 
with the following boundary conditions: 
$\hvec{q}=\hvec{0}$ on the substrate $y=0$\,; the normal stress on the free 
surface must balance normal surface tension, that is,
\begin{equation}
p=\hvec{\tilde{\tau}_{n}}\cdot\hvec{\tilde{n}}-\gamma\tilde{\kappa}\qquad
\mbox{on $y=\eta$},
\label{normstress}
\end{equation}
where $\hvec{\tilde{\tau}_{n}}$ is the deviatoric stress across the surface,
 $\tilde{\kappa}$ is the mean curvature of the free surface, 
$\hvec{\tilde{n}}$ is a unit normal to the free surface, $\gamma$ is the local
value of the surface tension and a tilde indicates evaluation at the free 
surface; tangential stress on the free surface must equal the surface tension 
gradients, that is,
\begin{equation}
\hvec{\tilde{\tau}_{n}}\cdot\hvec{\tilde{t}}=\hvec{\tilde{t}}\cdot\hvec{\nabla}
\gamma\qquad\mbox{on $y=\eta$},
\label{tangstress}
\end{equation}
where $\hvec{\tilde{t}}$ is a unit tangent to the free surface; and the 
kinematic condition 
\begin{equation}
\frac{\partial\eta}{\partial t}=\tilde{v}-\tilde{u}\frac{\partial\eta}
{\partial x}
\label{kinematic}
\end{equation}
states that the fluid particles on the free surface must 
follow the free surface.

The dynamics of a surfactant on the fluid surface is described by a PDE which
we derive here using conservation arguments. Consider an arbitrary interval of
substrate, $x$ in $I=[a,b]$, and the fluid above it. Let the fluid surface 
have a concentration of surfactant per unit area of the fluid surface, 
$\Gamma_{0}\Gamma(x,t)$ where $\Gamma_{0}$ is a typical value for the
concentration and $\Gamma(x,t)$ gives the nondimensional variations. Thus the 
nondimensional surface concentration per unit area of substrate is 
$\sqrt{1+\eta^{2}_{x}}\,\Gamma$. Conservation of mass implies 
that the rate of change of surfactant mass in $I$ is equal to the rate of 
mass influx across the ends. The rate of change of mass of surfactant in $I$ is
\begin{eqnarray}
\frac{d}{dt}\int_{I}\left[\Gamma\sqrt{1+\eta^{2}_{x}}\right]\,dx &=& \int_{I}
\left[\sqrt{1+\eta^{2}_{x}}\frac{\partial\Gamma}{\partial t}+\Gamma\frac{
\partial}{\partial t}\sqrt{1+\eta^{2}_{x}}\right]\,dx \\
\label{surf1}
&=& \int_{I}\left[\sqrt{1+\eta^{2}_{x}}\frac{\partial\Gamma}{\partial t}+
\Gamma\frac{\eta_{x}\eta_{xt}}{\sqrt{1+\eta^{2}_{x}}}\right]\,dx\,.
\label{surf2}
\end{eqnarray}
The fluid on the free surface, lateral velocity $\tilde{u}$, carries the
surfactant to produce a flux parallel to the substrate (in the
$x$-direction) of
\begin{equation}
\tilde{u}\Gamma{\sqrt{1+\eta^{2}_{x}}}\,.
\label{surf3}
\end{equation}
Molecular diffusion of surfactant on the free surface carries a flux
\begin{equation}
\frac{D_{s}}{\sqrt{1+\eta^{2}_{x}}}\frac{\partial\Gamma}{\partial s}=
\frac{D_{s}}{{(1+\eta^{2}_{x})}}\frac{\partial\Gamma}{\partial x}
\label{surf4}
\end{equation}
where $D_{s}$ is a surface diffusivity coefficient. The additional geometric
term, $1/\sqrt{1+\eta^{2}_{x}}$, appears because although the flux on the 
surface is $D_{s}\frac{\partial\Gamma}{\partial s}$ the direction of the flux 
is at an angle to the substrate and requires multiplication by the direction 
cosine $1/\sqrt{1+\eta^{2}_{x}}$\,.
The net rate of gain of surfactant in the interval $I$ through these transport
mechanisms is
\begin{equation}\small
\left.\left[\Gamma\tilde{u}{\sqrt{1+\eta^{2}_{x}}}-\frac{D_{s}}{{(
1+\eta^{2}_{x})}}\frac{\partial\Gamma}{\partial x}\right]\right|^{b}_{a} 
= \int_{I}\frac{\partial}{\partial x}\left[\Gamma\tilde{u}{\sqrt{1+
\eta^{2}_{x}}-\frac{D_{s}}{(1+\eta^{2}_{x})}\frac{\partial\Gamma}{\partial x}}
\right]\,dx\,.
\label{surf5}
\end{equation}
Equating (\ref{surf2}) to (\ref{surf5}) using (\ref{kinematic})
and after some algebraic manipulation, we obtain (noting $\tilde{u}_{x}$ is
$\frac{\partial}{\partial x}\left(\tilde{u}\right)$, not $\frac{\partial
u}{\partial x}\left|_{y=\eta}\right.$, and similarly for $\tilde{v}_{x}$)
{\begin{equation}
\frac{\partial\Gamma}{\partial
t}=\frac{\ep}{\sqrt{1+\eta_{x}^{2}}}
\frac{\partial}{\partial
x}\left[\frac{\Gamma_{x}}{(1+\eta_{x}^{2})}
\right]-\frac{\partial}{\partial x}\left(
\Gamma\tilde{u}\right)
+\frac{\Gamma\tilde{u}_{x}
\eta_{x}^{2}-\Gamma\tilde{v}_{x}\eta_{x}}{(1+\eta_{x}^{2})}\, ,
\label{surf6}
\end{equation}
where $\ep=1/\pe=D_{s}\mu/\gamma_{0} H$ is an inverse P\'{e}cl\'{e}t
number characterising the importance of surface diffusion compared
with advection.

To obtain a well-posed problem a further equation is needed relating
surface tension, $\gamma$, and 
surfactant concentration, $\Gamma$. It is well established~\cite[e.g.]
{adamson67,levich62,sheludko67} that the surface tension is a function of the 
surface concentration, $\gamma=\gamma(\Gamma).$ Here we chose to use 
a linear relationship between surface tension and surfactant concentration 
following Schwartz {\it et al}~\cite{schwartz95} namely,
\begin{equation} 
\gamma = 1+A(1-\Gamma)\,,
\label{surftenrel}
\end{equation}
where $\gamma$ has been nondimensionalised by scaling with respect to the
reference value of the surface tension $\gamma_{0}$ at $\Gamma=\Gamma_{0}$ and
\begin{equation}
   A=\frac{\Gamma_{0}}{\gamma_{0}}\frac{\partial\gamma}{\partial\Gamma}\,.
\label{constant}
\end{equation}

\section{The basis of the centre manifold analysis}
In this section we lay the basis for forming an accurate model of the 
dynamics of the thin fluid film with surfactant.
We adapt the governing fluid equations (\ref{continuity}--\ref
{navier}), the surfactant evolution equation (\ref{surf6}), the relevant 
boundary conditions (\ref{normstress}--\ref{kinematic}) and the constituent
equation (\ref{surftenrel}) to a form suitable for the application of centre
manifold theory and techniques in order to generate a low-dimensional 
model for the dynamics.

We develop a model of slow large scale flow by invoking the slowly varying
assumption, that is $\partial/\partial x$ is small, and with weak forcing, that
is $\bond$ and $\ham$ are also small. In centre manifold theory 
this is achieved by treating $\partial/\partial x$, $\mathcal{B}$ and 
$\mathcal{H}$ terms as ``nonlinear'' perturbations.
Thus the linear picture is obtained by neglecting any $\partial/\partial x$,
$\mathcal{B}$ and $\mathcal{H}$ terms. This may be seen as being equivalent to 
the multiple-scale assumption of variations occurring only on a large lateral 
length scale (see Roberts~\cite{roberts88,roberts96,roberts97} for a fuller 
explanation). We also assume that the fluid is thin enough for
gravity to be a perturbing influence but thick enough for van der Waals forces
to also be a perturbing influence, that is, the Bond and Hamaker numbers are
both small.

The ``linear'' dynamics are then solutions of the following equations
\begin{eqnarray}
   \frac{\partial v}{\partial y} &=& 0\,, \\
   \label{lincont}
   \mathcal{R}\frac{\partial\hvec{q}}{\partial t}+\frac{\partial p}{\partial y}
   \hvec{j}-\frac{\partial^{2}\hvec{q}}{\partial y^{2}} &=& 0\,, \\
   \label{linns}
   \frac{\partial\Gamma}{\partial t} &=& 0\,,
\end{eqnarray}
with boundary conditions
\begin{eqnarray}
   \hvec{q} &=& \hvec{0}\quad\mbox{on $y=0$}\,, \\
   -p+2\frac{\partial v}{\partial y} &=& 0\quad\mbox{on $y=\eta$}\,, \\
   \frac{\partial u}{\partial y}+\frac{\partial v}{\partial x} &=& 0\quad\mbox
   {on $y=\eta$}\,, \\
   \frac{\partial\eta}{\partial t}-v &=& 0\quad\mbox{on $y=\eta$}\,,
\end{eqnarray}
All solutions of these linear equations are composed of the decaying
lateral shear modes $v=p=0,u=b\sin{(l\pi y/(2\eta))} \exp{(\lambda_{l} t)}$, 
together with two critical modes $\eta=$  constant and $\Gamma=$ constant. 
Here the integer $l$ parameterises the vertical wavenumber and the
decay rate of the lateral shear modes are
$-\lambda_{l}=l^{2}\pi^{2}/(4\eta^{2}\re)$. So linearly, and in
the absence of any lateral variations on a flat substrate, all the lateral
shear modes decay exponentially quickly, on a time-scale of
$\re\eta^{2}$, just leaving a film of constant thickness with
a covering of surfactant of constant concentration as the permanent mode. This
spectrum, of all eigenvalues strictly negative except for a few
that are zero, is the classic spectrum for the application of centre
manifold theory: the Existence Theorem in~\cite{carr81} assures that
the nonlinear effects in the physical equations just perturb this linear 
picture of the dynamics so that in the long-term all solutions of the full 
nonlinear system are dominated by the slow dynamics induced by nonlinearities
and large-scale lateral
variations in the film thickness and contaminant. The Relevance Theorem
in~\cite{carr81} assures that these dynamics are exponentially
attractive, asymptotically complete, and so form a generic model of
the long-term dynamics of the contaminated film. With the caveat that a
strict theory has not yet been developed to cover this application to 
non-linear large-scale flows,
the closest being that of Gallay~\cite{gallay93} and H\u{a}r\u{a}gus~
\cite{haragus95} (but also see~\cite{roberts88}), the centre manifold 
concepts and techniques are applied to systematically develop a 
low-dimensional lubrication model of the dynamics of the film.

Having identified the critical modes associated with the zero
decay-rate, the subsequent analysis is straightforward. The usual
approach is to write the fluid fields $\hvec{q}(x,y,t)=(u,v)$ and $p(x,y,t)$,
as a function of the critical modes $\eta$ and $\Gamma$ (equivalent to the
``slaving'' principle of synergetics~\cite{haken93}). Instead of seeking
\emph{explicit} asymptotic expansions in the ``amplitudes'' of the critical 
modes~\cite{roberts88,roberts96}, an iterative algorithm is applied to find 
the centre manifold and the evolution thereon which is based directly upon the
Approximation Theorem in~\cite{carr81,roberts97} and its variants;
explained in detail by Roberts in~\cite{roberts97}.

\section{The centre manifold model}
We solve the continuity and Navier-Stokes equations under the assumptions 
introduced above by programming the computer algebra package \textsc{reduce}. 
The program listed in the Appendix iteratively solves the physical equations
using techniques explained by Roberts~\cite{roberts88,roberts97}.
We express the velocity and pressure fields in terms of the scaled normal 
coordinate $\zeta=y/\eta(x,t)$ to simplify the expressions; in this stretched 
coordinate the free surface is $\zeta=1$. The computer algebra gives the fluid
fields to be
\begin{eqnarray}
   u &\approx& \left(\zeta-\frac{1}{2}\zeta\right)\bond\sin{\theta}\,\eta^{2}+
   \left(3\zeta-\frac{3}{2}\zeta^{2}\right)\ham\frac{\eta_{x}}
   {\eta^{2}}\nonumber \\
   && {}-\left(\zeta-\frac{1}{2}\zeta^{2}\right)\bond\cos{\theta}\,
   \eta^{2}\eta_{x}+\left(\frac{5}{2}\zeta-\frac{1}{2}\zeta^{2}-\frac{1}{3}
   \zeta^{3}\right)\bond\sin{\theta}\,\eta^{3}\eta_{xx}\nonumber \\
   && {}+\left(\zeta-\frac{1}{2}\zeta^{2}\right)\gamma\eta^{2}\eta_{xxx}
   +\zeta\eta\gamma_{x}\,, \\
   v &\approx& -\frac{1}{2}\zeta^{2}\bond\sin{\theta}\,\eta^{2}\eta_{x}+\left(
   \frac{9}{2}\zeta^{2}-2\zeta^{3}\right)\ham\frac{\eta^{2}_{x}}{\eta^{2}}+
   \frac{1}{2}\zeta^{2}\bond\cos{\theta}\,\eta^{2}\eta_{x}^{2}\nonumber \\
   && {}-\left(\frac{3}{2}\zeta^{2}-\frac{1}{2}\zeta^{3}\right)\ham
   \frac{\eta_{xx}}{\eta}+\left(\frac{1}{2}\zeta^{2}-\frac{1}{6}\zeta^{3}
   \right)\bond\cos{\theta}\,\eta^{3}\eta_{xx}\nonumber \\
   && {}-\frac{1}{2}\zeta^{2}\eta^{2}\gamma_{xx}-\frac{5}{2}\zeta^{2}\bond
   \sin{\theta}\,\eta^{2}\eta^{3}_{x}-\left(\frac{15}{2}\zeta^{2}-\frac{1}{2}
   \zeta^{3}\right)\bond\sin{\theta}\,\eta^{3}\eta_{x}\eta_{xx}\nonumber \\
   && {}-\left(\frac{5}{4}\zeta^{2}-\frac{1}{6}\zeta^{3}-\frac{1}{12}\zeta^{4}
   \right)\bond\sin{\theta}\,\eta^{4}\eta_{xxx}\,, \\
   p &\approx& -\gamma\eta_{xx}+(1-\zeta)\bond\cos{\theta}\,\eta-\left(1+\zeta
   \right)\bond\sin{\theta}\,\eta\eta_{x}+(1+\zeta)\bond\cos{\theta}\,\eta
   \eta^{2}_{x}\nonumber \\
   && {}+\left(\frac{1}{2}+\zeta-\frac{1}{2}\zeta^{2}\right)\bond\cos{\theta}
   \,\eta^{2}\eta_{xx}+\left(3+9\zeta-6\zeta^{2}\right)\ham\frac{\eta_{x}^{2}}
   {\eta^{3}}\nonumber \\
   && {}-\left(\frac{3}{2}+3\zeta-\frac{3}{2}\zeta^{2}\right)\ham\frac{
   \eta_{xx}}{\eta^{2}}-\left(9+5\zeta\right)\bond\sin{\theta}\,\eta
   \eta^{3}_{x}\nonumber \\
   && {}-\left(\frac{27}{2}+15\zeta-\frac{3}{2}\zeta^{2}\right)\bond\sin
   {\theta}\,\eta^{2}\eta_{x}\eta_{xx}\nonumber \\
   && {}-(1-\Gamma)A\eta_{xx}-2\eta_{x}\gamma_{x}-\left(1+\zeta\right)\eta
   \gamma_{xx}\,.
\end{eqnarray}
See that the lateral velocity, $u$, is approximately parabolic, Poiseuille 
flow, which form components of the forcing that act through lateral pressure 
gradients, but is linear, Couette flow, from the surface tension gradients. 
Then the velocity normal to the substrate, $v$, follows predominately from the
continuity equation. These expressions give comprehensive details of the 
physical fields corresponding to any particular $\eta(x,t)$ and $\Gamma(x,t)$.

The computer algebra also derives the corresponding evolution equations for 
$\eta$ and $\Gamma$ which are
\begin{eqnarray}
    \frac{\partial\eta}{\partial t} & = & -\frac{1}{2}\frac{\partial}
    {\partial x}\left(\eta^{2}\gamma_{x}\right)-\frac{1}{3}\frac{\partial}
    {\partial x}\left(\eta^{3}\left(\gamma\eta_{xx}\right)_{x}\right)
    \nonumber \\
    & & {}-\frac{\partial}{\partial x}\left(\frac{1}{3}\eta^{3}+\frac{7}{3}
    \eta^{3}\eta^{2}_{x}+\eta^{4}\eta_{xx}\right)\bond\sin{\theta}\nonumber \\
    & & {}+\frac{\partial}{\partial x}\left(\frac{1}{3}\eta^{3}\eta_{x}
    +\frac{3}{5}\eta^{5}\eta_{xxx}+4\eta^{4}\eta_{x}\eta_{xx}+\frac{7}{3}
    \eta^{3}\eta_{x}^{3}\right)\bond\cos{\theta}\nonumber \\
    & & {}+\frac{\partial}{\partial x}\left(-\frac{\eta_{x}}{\eta}+
    \frac{48}{5}\eta_{x}\eta_{xx}-\frac{9}{5}\eta\eta_{xxx}-7\frac{
    \eta_{x}^{3}}{\eta}\right)\ham\nonumber \\
    & & {}+\frac{\partial}{\partial x}\left(\frac{32}{105}\eta^{2}\eta^{2}_{x}
    -\frac{10}{21}\eta^{3}\eta_{xx}\right)\ham\re\bond\sin{\theta}\nonumber \\
    & & {}+\frac{\partial}{\partial x}\left(\frac{44}{105}\eta^{3}
    \eta_{x}\eta_{xx}+\frac{4}{15}\eta^{4}\eta_{xxx}-\frac{4}{105}\eta^{2}
    \eta_{x}^{3}\right)\ham\re\bond\cos{\theta}\nonumber \\
    & & {}+\mathcal{O}(\partial^{6}_{x},\bond^{2},\ham^{2})\quad\mbox{and}
    \label{eta} \\
    \frac{\partial\Gamma}{\partial t} & = & -\frac{\partial}{\partial x}
    \left(\Gamma\eta\gamma_{x}\right)-\frac{1}{2}\frac{\partial}{\partial x}
    \left(\Gamma\eta^{2}\left(\gamma\eta_{xx}\right)_{x}\right)\nonumber \\
    && {}+\frac{\partial}{\partial x}\left(-\frac{1}{2}\Gamma\eta^{2}-
    \frac{5}{3}\Gamma\eta^{3}\eta_{xx}-\frac{17}{4}\Gamma\eta^{2}\eta_{x}^{2}
    \right)\bond\sin{\theta}\nonumber \\
    && {}+\left(\frac{3}{2}\Gamma\eta\eta_{x}^{3}-\frac{1}{4}\Gamma_{x}
    \eta^{2}\eta_{x}^{2}\right)\bond\sin{\theta}\nonumber \\
    && {}+\frac{\partial}{\partial x}\left(\frac{1}{2}\Gamma\eta^{2}\eta_{x}+
    4\Gamma\eta^{2}\eta_{x}^{3}+\frac{20}{3}\Gamma\eta^{3}\eta_{x}\eta_{xx}+
    \Gamma\eta^{4}\eta_{xxx}\right)\bond\cos{\theta}\nonumber \\
    && {}+\left(-\Gamma\eta\eta_{x}^{4}+\frac{1}{3}\Gamma\eta^{3}\eta_{xx}^{2}+
    \frac{1}{2}\Gamma_{x}\eta^{2}\eta^{3}_{x}+\frac{1}{3}\Gamma_{x}\eta^{3}
    \eta_{x}\eta_{xx}\right)\bond\cos{\theta}\nonumber \\
    && {}+\frac{\partial}{\partial x}\left(-\frac{3}{2}\frac{\Gamma
    \eta_{x}}{\eta^{2}}-\frac{32}{3}\frac{\Gamma\eta_{x}^{3}}{\eta^{2}}+16\frac
    {\Gamma\eta_{x}\eta_{xx}}{\eta}-3\Gamma\eta_{xxx}\right)\ham\nonumber \\
    && {}+\left(-\frac{1}{3}\frac{\Gamma\eta_{x}^{4}}{\eta^{3}}-\frac{\Gamma
    \eta^{2}_{xx}}{\eta}+\frac{7}{6}\frac{\Gamma_{x}\eta_{x}^{3}}{\eta^{2}}
    -\frac{\Gamma_{x}\eta_{x}\eta_{xx}}{\eta}\right)\ham\nonumber \\
    & & {}+\frac{\ep}{\sqrt{1+\eta_{x}^{2}}}\frac{\partial}{\partial x}\left(
    \frac{\Gamma_{x}}{1+\eta_{x}^{2}}\right)\nonumber \\
    && {}+\frac{\partial}{\partial x}\left(-\frac{89}{120}\Gamma\eta^{2}
    \eta_{xx}+\frac{7}{15}\Gamma\eta\eta_{x}^{2}\right)\ham\re\bond
    \sin{\theta}\nonumber \\
    && {}+\frac{\partial}{\partial x}\left(\frac{13}{20}\Gamma
    \eta^{2}\eta_{x}\eta_{xx}+\frac{5}{12}\Gamma\eta^{3}\eta_{xxx}-\frac{1}{20}
    \Gamma\eta\eta_{x}^{3}\right)\ham\re\bond\cos{\theta}\nonumber \\
    && {}+\mathcal{O}(\partial^{6}_{x},\bond^{2},\ham^{2}).
    \label{gamma}
\end{eqnarray}
The error term $\mathcal{O}(\partial^{6}_{x},\bond^{2},\ham^{2})$ indicates the
terms retained in the model by neglecting any term with 6 or more spatial
derivatives or any quadratic or higher terms in $\bond$ or $\ham$. Thus we
retain the terms seen above. Any particular application need not retain all of
these terms. The Approximation Theorem in \cite{carr81} supports many 
consistent truncations of
these expressions and the retention of terms is only dependent on the type of
application. This model is comprehensively flexible in that it encompasses any
model of a similar genre such as those described below.

Gaver and Grotberg~\cite{gaver90}, de Wit {\it et al}~\cite{dewit94}, 
Schwartz and Weidner~\cite{schwartz95} and others have previously developed 
evolution models for this flow using heuristic arguments based on traditional 
lubrication theory. These models are all similar to the de Wit model:
\begin{eqnarray}
    \frac{\partial\eta}{\partial t} & \approx & -\frac{1}{2}\frac{\partial}
    {\partial x}\left(\eta^{2}\gamma_{x}\right)-\frac{1}{3}\frac{\partial}
    {\partial x}\left(\eta^{3}\left(\gamma\eta_{xx}\right)_{x}\right)
    -\ham\frac{\partial}{\partial x}\left(\frac{\eta_{x}}{\eta}\right)
    \quad\mbox{and}
    \label{dewiteta} \\
    \frac{\partial\Gamma}{\partial t} & \approx & -\frac{\partial}{\partial x}
    \left(\Gamma\gamma_{x}\eta\right)-\frac{1}{2}\frac{\partial}{\partial x}
    \left(\Gamma\eta^{2}\left(\gamma\eta_{xx}\right)_{x}\right)
    -\frac{3\ham}{2}\frac{\partial}{\partial x}\left(\frac{\Gamma\eta_{x}}
    {\eta^{2}}\right)+\ep\Gamma_{xx},
    \label{dewitgamma}
\end{eqnarray}
and we consider them to be a subset of ours. We have included the effects of
gravity (indicated by the Bond number); steeper surface slopes (through
$\footnotesize{\sqrt{1+\eta_{x}^{2}}}$) and the interaction between gravity 
and van der Waals forces (indicated by $\bond\ham$) which also involves 
inertia as the Reynolds number appears.
Numerical solutions are compared in $\S 6$.

\section{Stability analysis of simple flows}
Linearising the two model's evolution equations~(\ref{eta}--\ref{gamma}) about
a fixed point gives insight into the physical effects taking place. 
Assume the fluid film is flat with the thickness of 
the film and the average surfactant concentration both being one in 
nondimensional units. Then we elucidate the interactive 
dynamics of the system by perturbing these values. We assume a 
solution of these equations with initial sinusoidal ripples with growth rate 
$\lambda$ or decay rate $-\lambda$ and lateral wavelength $k$ of the form
\begin{equation}
    \eta = 1+ae^{(\lambda t+ikx)}\quad\mbox{and}\quad
    \Gamma = 1+be^{(\lambda t+ikx)}
    \label{lineqns}
\end{equation}
for some $a,b$.

\begin{figure}
\begin{center}\includegraphics[scale=0.7]{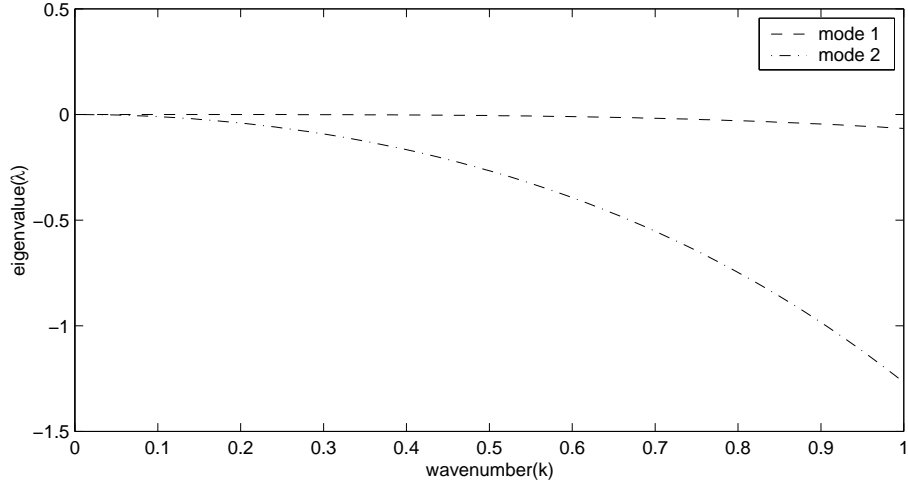}\end{center}
\caption{The eigenvalues of the linear modes~(\ref{lineqns}) are
plotted against the wavenumber $k$. This shows one mode decaying rapidly and
one mode decaying very slowly. This slowly decaying mode represents physically
long lasting corrugations on the surface maintained by in phase surfactant 
variations. The inverse P\'{e}cl\'{e}t number ($\ep=1/\pe$) chosen here is the
typical value $10^{-4}$.}
\label{modefig}
\end{figure}

Substitute these expressions into the model (\ref{eta}--\ref{gamma}) and 
neglect all nonlinear terms in $a$ and $b$ to produce
\begin{eqnarray}
    \left(\lambda+\frac{k^{4}}{3}\right)a & = & -\left(
    \frac{k^{2}}{2}\right)b\,, \\
    \label{thirteen}
    \left(\lambda+\ep k^{2}+k^{2}\right)b & = & 
    -\left(\frac{k^{4}}{2}\right)a\,.
    \label{fourteen}
\end{eqnarray}
Nontrivial solutions only exist when
\begin{equation}
\lambda^{2}+\left(k^{2}+\ep k^{2}+\frac{1}{3}k^{4}\right)\lambda+\frac{1}{12}
k^{6}+\frac{1}{3}\ep k^{6}=0\,.
\label{chareqn}
\end{equation}
The eigenvalues of the linear modes~(\ref{lineqns}),
\begin{eqnarray}
   \lambda &=& -\left(\frac{1}{2}+\frac{1}{2}\ep+\frac{1}{6}k^{2}\right)k^{2}
   \nonumber \\
   &&{}\pm\sqrt{\frac{1}{4}k^{4}+\frac{1}{2}\ep k^{4}+\frac{1}{12}k^{6}+
   \frac{1}{4}\ep^{2}k^{4}-\frac{1}{6}\ep k^{6}+\frac{1}{36}k^{8}}\,,
   \label{eigen}
\end{eqnarray}
are plotted in Figure~\ref{modefig}. They show that one mode decays rapidly
with respect to the other which decays very slowly. This last mode denotes 
physically long lasting corrugations on the surface maintained by in phase 
surfactant variations.

This linear analysis of the dynamical system shows that the model is stable 
for all wave numbers $k$. Therefore the model (\ref{eta}--\ref{gamma}) is 
structurally stable.

\section{Numerical simulations}

\begin{figure}
\begin{center}\includegraphics{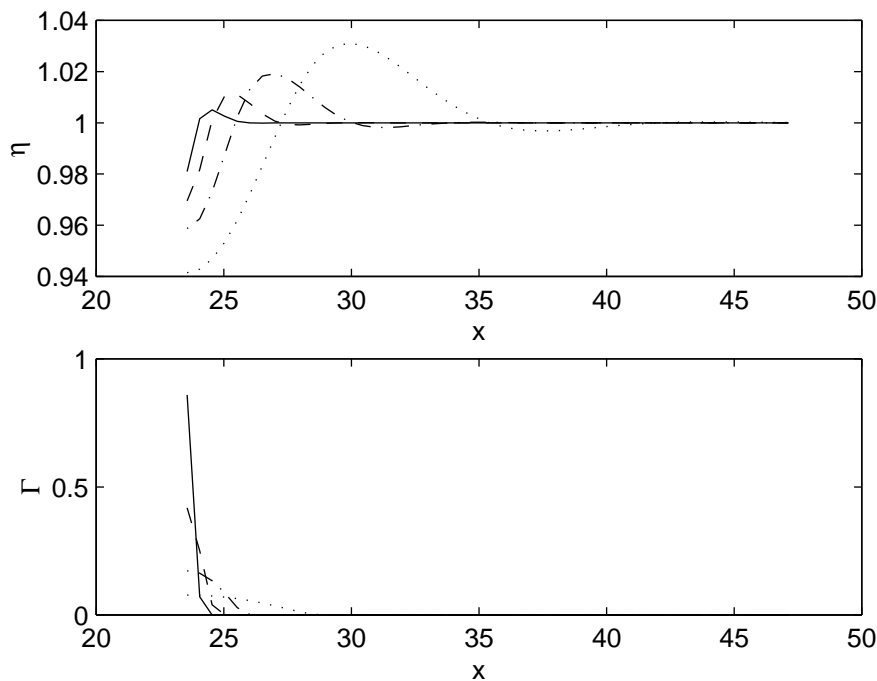}\end{center}
\caption{The film thickness (top) and surfactant concentration
(bottom) profiles at non--dimensional times of $t=1$ (---),
$t=10$ (- -), $t=100$ ($-\cdot$) and $t=1000$ ($\cdots$). The initial
conditions are a flat free surface with a drop of surfactant placed
at $x=15\pi/2$. Symmetry is assumed.}
\label{flatevol}
\end{figure}

Our model (\ref{eta}--\ref{gamma}) and the de Wit model 
(\ref{dewiteta}--\ref{dewitgamma}) are solved for $\eta(x,t)$ and 
$\Gamma(x,t)$ using a standard $1^{st}$ order backward time and $2^{nd}$ order
centred spatial differencing~\cite{degrez92}. The properties of the fluid are 
taken to be surface tension $\gamma=30\, \mbox{dynes/cm}$, viscosity 
$\mu=10^{-2}\, \mbox{g/(cm\,s)}$, density $\rho=1\, \mbox{g/}\mbox{cm}^{3}$ 
and the surface diffusivity constant $D_{s}= 10^{-4}\, \mbox{cm}^{2}\mbox{/s}$
with the uniform thickness of the film $10^{-5}$\, cm and a drop of surfactant
of $10^{-10}\mbox{mol/cm$^2$}$ placed initially in the centre of the fluid 
surface. This corresponds to a Reynolds number of $\re=3$, a Bond number of
$\bond=3\times10^{-11}$, a Hamaker constant of $\ham=0.001$ and a 
P\'{e}cl\'{e}t number of $\pe=1/\ep=300$. 
Therefore it is expected that the flow will be dominated by viscous forces and 
that the transport of surfactant will be dominated by advection.

\begin{figure}
\begin{center}\includegraphics[scale=0.7]{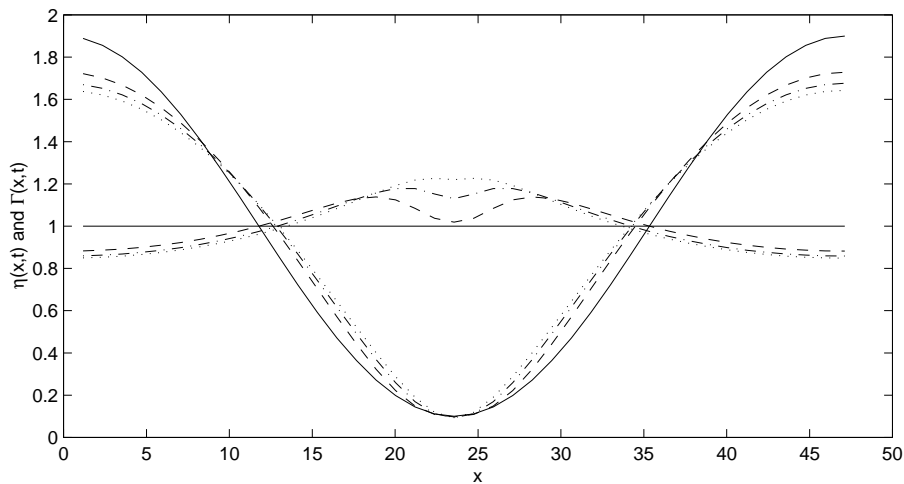}\end{center}
\caption{The film thickness (top) and surfactant concentration (bottom)
profiles at non-dimensional times of $t=0$ (---), $t=15$ (- -), $t=30$ 
($-\cdot$) and $t=45$ ($\cdots$). The initial conditions are a corrugated free 
surface with a even layer of surfactant on the fluid.}
\label{smallcorr}
\end{figure}

For simualation over a large spatial domain, convergence was obtained with a
single Newton step with a non-dimensional time step of $\delta t=100$ was on a
spatial grid of $N=97$ points. This numerical scheme is stable for all 
time-steps (except when ludicrously big) and spatial grids. The time 
step of 100 is small enough so that the dynamics of a reasonable 
number of spatial modes are modelled accurately by the scheme.

Figure \ref{flatevol} shows the temporal evolution of the film thickness
and surfactant concentration for our centre manifold
model (\ref{eta}--\ref{gamma}). The first graph
plots $\eta$, the film thickness, as a function of $x$ for various times
whilst the second graph plots $\Gamma$, the surfactant concentration, also as a
function of $x$ for various times. Surface tension gradients formed
due to the action of the surfactant result in the propagation of a
front. The initially steep concentration gradients die out over time
due to the advection of the surfactant by the front.

The stability analysis in $\S$5 shows that there could be long lasting
corrugations on the surface maintained by in phase surfactant variations.
Figures~\ref{smallcorr} and~\ref{corrugation} show plots of the 
evolution of a fluid with an initially corrugated free surface contaminated 
with a even layer of surfactant. The initial deformation in Figure 
\ref{smallcorr} shows the accumulation of surfactant in the trough caused by
the pressure gradients driving fluid into the trough. After a non-dimensional 
time step of $t=45$ the raised 
concentration of surfactant generates surface tension gradients to oppose this 
collapse and leads to the corrugations lasting a long time as shown in Figure
\ref{smallcorr}. This agrees with the stability analysis in Section 5 and 
demonstrates that surfactants hinder the levelling of thin films.

\begin{figure}
\begin{center}\includegraphics[scale=0.7]{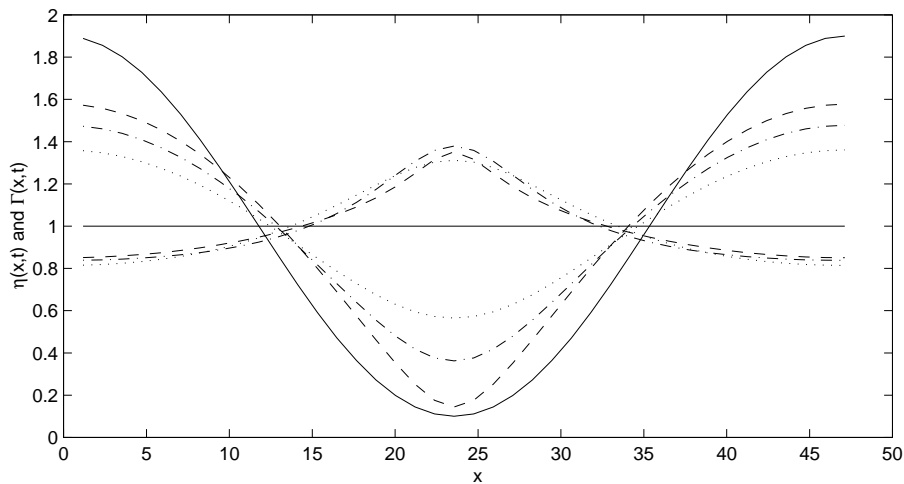}\end{center}
\caption{The film thickness (top) and surfactant concentration (bottom)
profiles at non-dimensional times of $t=0$ (---), $t=100$ (- -), $t=200$ 
($-\cdot$) and $t=300$ ($\cdots$). The initial conditions are a corrugated
free surface with a even layer of surfactant on the fluid.}
\label{corrugation}
\end{figure}

Our model for the temporal surfactant evolution (\ref{gamma}) when compared
with the de Wit model, contains additional fourth order terms involving 
$\eta_{x}^{2}$. The enhanced accuracy of the model
becomes apparent when the surface gradients, $\eta_{x}$, are
sufficiently large. Figure \ref{comparison} shows a comparison between
our model and the de Wit model. The difference between the two
solutions is plotted at a common nondimensional time of $t=10$ and
shows the increasing disparity as the P\'{e}cl\'{e}t number is decreased.

\begin{figure}
\begin{center}\includegraphics[scale=0.7]{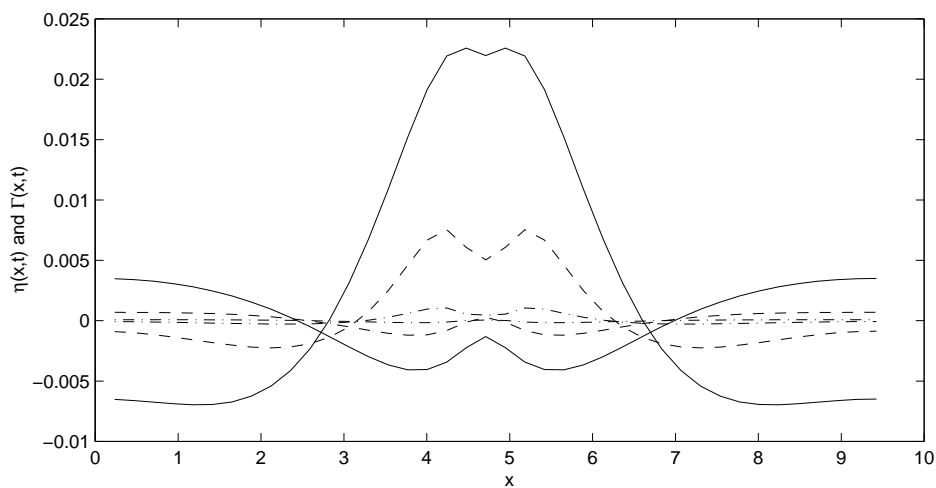}\end{center}
\caption{The film thickness (top) and surfactant concentration
(bottom) profiles at the same nondimensional time, $t=10$. The
graphs plot the difference between the centre manifold model,
(\ref{eta}--\ref{gamma}), and the de Wit model (\ref{dewiteta}--\ref
{dewitgamma}), at different P\'{e}cl\'{e}t numbers $\pe=3$ (---), 
$\pe=30$ (- -) and $\pe=300$ ($-\cdot$). 
As the P\'{e}cl\'{e}t numbers decrease the difference between 
the models increase.}
\label{comparison}
\end{figure}

\section{Summary}
In this paper we have developed a comprehensive structurally stable model of 
the dynamics of the spreading of a contaminant on a thin fluid. The centre 
manifold approach we adopt incorporates all the physical effects at the 
appropriate stage in the modelling. Numerical simulations have
demonstrated the importance of the extra terms in our model compared
with other models developed using other methods.  

\bibliography{fluids}
\bibliographystyle{plain}

\appendix
\section{Computer algebra code}
A computer algebra program to perform all the necessary detailed
algebra for this physical problem was used. An important feature
of this iteration is that it is performed until the residuals of the
actual governing equations are zero, to some order of error. Thus the
correctness of the results is based only upon the correct evaluation
of the residuals and sufficient iterations.

\subsection*{}

\listinginput{1}{program.red}

\newpage

\end{document}